\documentclass[]{article} 
\usepackage[left=2.5cm,centering]{geometry}
\providecommand{\keywords}[1]
{
  \small	
  \textbf{\textit{Keywords -- }} #1
}
\usepackage[english]{babel}
\usepackage{mathtools}
\usepackage{mathrsfs}
\usepackage{cancel}
\usepackage{graphicx}
\usepackage{empheq}
\usepackage{amsfonts,amssymb,amsmath}
\usepackage{enumitem}
\usepackage{authblk}
\usepackage[table]{xcolor}
\usepackage{tikz}
\usepackage{pgfplots}
\usepackage{bm}    
\usepackage{diagbox}
\usepackage{booktabs}

\usepackage{natbib}
\usepackage{algorithm}

\newtheorem{theorem}{Theorem}
\newcommand{\al}{$\alpha\,$}

\begin{document}

\title{A Neural-Network Approach for the Stability of Loosely-Coupled FSI Simulations in Hemodynamics.}

\author[1,3]{Francesca Renzi\thanks{ corresponding author}}
\author[2]{Stefano Pagani}
\author[2]{Paola F. Antonietti}
\author[1]{Christian Vergara}

\affil[1]{LaBS, Department of Chemistry, Materials, and Chemical Engineering, Politecnico di Milano, Milano, Italy\\
\texttt{francesca.renzi@polimi.it, christian.vergara@polimi.it}}
\affil[2]{MOX, Department of Mathematics, Politecnico di Milano, Milano, Italy\\
\texttt{stefano.pagani@polimi.it, paola.antonietti@polimi.it}}
\affil[3]{Inria, Paris, France}
\date{}

\maketitle

\begin{abstract}
The development of accurate and efficient computational methods for fluid-structure interaction (FSI) problems remains an active area of research, given the wide range of engineering and biomedical applications involving FSI.
Among the available approaches, partitioned loosely coupled (LC) schemes are particularly attractive for large-scale simulations due to their modularity and reduced computational cost. 
However, their applicability is limited by stability issues in high added-mass regimes, where the fluid and structure densities are comparable, as occurs in hemodynamics.
In our recent work, we proposed an LC algorithm based on Dirichlet-Neumann interface coupling, referred to as LC-DN-\al scheme, which improved the stability of the coupling by introducing a stabilization parameter \al. 
Although theoretical stability has been demonstrated under suitable constraints on \al, the dependence of the optimal value of \al on the physical and numerical characteristics of general FSI problems is not yet fully understood, making its practical selection challenging.

In this work, we address this open question by proposing a data-driven strategy for automatically estimating the optimal value of \al.
Specifically, we develop a feed-forward neural network to predict \al from the descriptors of the FSI problem.
The model is trained on a dataset of baseline FSI problems in large added-mass regimes and is subsequently evaluated on more realistic configurations to assess its robustness and generalization capabilities.
The numerical results show that the proposed model generalizes across the considered FSI test scenarios and provides robust estimates of \al under moderate variations in the problem setup, including changes in geometry and loading conditions, that do not completely alter the fluid-structure interaction mechanics with respect to those represented in the training set. 

The methodology has been intentionally developed and validated on simplified FSI problems, providing a proof of concept for data-driven strategies to improve the efficiency of partitioned FSI simulations and supporting their extension to more complex and realistic scenarios.


\end{abstract}

\section{Introduction}

Fluid-structure interaction (FSI) is a common phenomenon in engineering, occurring across a wide range of applications.
These encompass both industrially relevant problems and topics of significant scientific interest, including aerodynamics, microfluidics, structural engineering, and biomechanics.
However, the numerical modeling of FSI problems remains demanding due to the nonlinear nature of the coupling between fluid and structural dynamics, the possible intrinsic nonlinearities within each physical subsystem, and the severe ill-conditioning of the linear systems arising from the discretization of the coupled problem \citep{richter2015monolithic}.
Consequently, the development of accurate and efficient numerical methods for FSI simulations is of key importance for addressing the computational challenges of real-world applications.
To this end, several strategies have been proposed and can be broadly classified into {\sl physics-based} and {\sl data-driven} methods \citep{tiba2025machine}.

Physics-based numerical methods continue to be the dominant approach for solving FSI problems as they accurately capture underlying physical phenomena and provide reliable predictions.
In this area, the modularity of {\sl partitioned} algorithms offers a significant simplification in the implementation complexity compared with {\sl monolithic} approaches \citep{farhat2000two, gee2011truly, crosetto2011parallel, richter2015monolithic, landajuela2017coupling, bucelli2022partitioned}. 
In particular, the possibility of solving the fluid and structure subproblems separately allows us to handle smaller, better-conditioned systems and reuse standalone physics-specific solvers and preconditioners, which is often required when using commercial software packages.

The fluid and structure solvers communicate at the fluid-structure interface through the enforcement of continuity conditions on relevant quantities (i.e., velocities and stresses). 
In particular, the classical {\sl Dirichlet-Neumann} (DN) interface coupling, which consists of imposing interface velocity as a Dirichlet boundary condition for the fluid problem and the interface stresses as a Neumann boundary condition for the structure problem, represents the most natural partitioned formulation, being also usually the only option implemented in commercial software.
Depending on how these conditions are enforced, the resulting algorithm is said {\sl strongly-coupled} (SC), through interface iterations \citep{farhat2000two, nobile2001numerical, matthies2002partitioned, causin2005added, badia2008fluid, kuttler2008fixed, degroote2010simulation, landajuela2017coupling, schussnig2022efficient}, or {\sl loosely-coupled} (LC), by performing a single fluid and structure solve per time step \citep{gigante2021stability,banks2014analysis,burman2014explicit,bukavc2013fluid,guidoboni2009stable}.

LC-DN schemes are particularly appealing for large-scale FSI simulations because they retain the modularity of partitioned methods, the simplicity of the DN coupling, and a substantial reduction in computational cost. 
Their main drawback is that they suffer from stability issues when the fluid and structure densities are comparable \citep{causin2005added, forster2007artificial, degroote2008stability, brummelen2009added}.
This phenomenon, known as the {\sl added-mass effect} and originating from the strong inertial interaction between the fluid and the structure, is especially pronounced in hemodynamics simulations, limiting the applicability of LC-DN schemes in this field. 

Extensive work has been done in the literature to improve the stability properties of LC algorithms in the context of large added-mass.
Most of the proposed approaches depart from the classical DN formulation by introducing alternative interface conditions, such as Robin- or Nitsche-type couplings, or by combining the partitioned strategy with operator-splitting techniques and suitable stabilization terms \citep{guidoboni2009stable, bukavc2013fluid, fernandez2013incremental, lukavcova2013kinematic,  banks2014analysis,  burman2014explicit, gigante2021stability, burman2025robin}.
Instead, in our previous work, we retained the classical DN coupling and proposed an LC algorithm ({\sl LC-DN-\al} scheme) whose stability is proved for suitable ranges of a stabilization parameter \al \citep{renzi2026stable, renzi2026choice}. 
However, estimating such ranges as well as the optimal value of \al for realistic simulations remains an open challenge.

In recent years, data-driven strategies based on Scientific Machine Learning (SciML) have been developed for the numerical treatment of FSI problems, thanks to advances in computational capabilities and the increasing availability of high-quality data.
Current SciML research in FSI has mainly focused on two directions: optimization problems \citep{shirvani2023machine} and the development of surrogate models, including reduced-order models and physics-informed neural networks, to approximate or replace high-fidelity simulations \citep{xiao2016non, shinde2019galerkin, dupont2023data, tiba2024non, zhang2025parametrized, farea2025learning, liu2026fluid}.
A third emerging line of research seeks to enhance the efficiency of conventional solvers rather than replacing them to reduce the cost of the coupling in the context of partitioned algorithms: by accelerating the convergence of fixed-point iterations \citep{totounferoush2021partitioned, tiba2025machine} or by improving the prediction of interface quantities \citep{mazhar2023novel}.
Despite promising, such strategies have so far been designed almost exclusively for SC schemes, and the potential of SciML approaches remains largely unexplored in the context of LC schemes. 

In this work, we introduce a feed-forward neural network model for the automatic estimation of the stabilization parameter \al that appears in the LC-DN-\al algorithm proposed in \cite{renzi2026choice, renzi2026stable}.
The network is trained on a set of baseline 2D FSI problems in large added-mass regimes, and the learned parameters are subsequently transferred to more realistic 2D configurations to evaluate their robustness and generalization beyond the training scenarios. 
To the best of our knowledge, this is the first study to integrate SciML techniques with LC algorithms for FSI simulations in high-added-mass regimes, such as those encountered in hemodynamics.
The study is intended as a proof of concept, demonstrating the potential of SciML for estimating the parameters of LC algorithms based solely on the physical and numerical characteristics of the problem at hand. 

\section{Algorithm formulation and stability of the LC-DN-\al scheme}
In this section, we review the LC FSI scheme proposed in \cite{renzi2026stable} together with the corresponding stability theoretical results.


\subsection{Formulation of the scheme}\label{subsec:formulation}
Let us consider the FSI problem resulting from the coupling between an incompressible Newtonian fluid and a linear elastic structure.

We formulate the fluid dynamics in an Arbitrary Lagrangian-Eulerian (ALE) framework on the time-dependent domain $\Omega_f^t$ and structure dynamics in a Lagrangian framework on its reference configuration $\Omega_s^0$. 
The governing equations consist of the incompressible Navier-Stokes (NS) equations in ALE form for the fluid unknowns, velocity $\bm u$ and pressure $p$, and the linear elasticity equations for the structure unknowns, displacement $\bm \eta$ and interface velocity $\bm w$, equipped with suitable initial and boundary conditions \citep{quarteroni2017cardiovascular}.
The two subproblems are coupled through the interface $\Sigma^t = \overline{\Omega}_f^t \cap \overline{\Omega}_s^t$ by enforcing the continuity of velocity and the equilibrium of interface tractions. 
Throughout the paper, we denote by $\bm{n}$ the unit normal vector to the interface, oriented outward with respect to the fluid domain.

The problem is discretized in time by means of the Finite-Difference method on the uniform partition $t^n = n\Delta t$, with $n = 0, 1, \dots, T/\Delta t$, of temporal domain $[0, T]$, where $\Delta t$ denotes the time step. 
In this work, we adopt an explicit treatment of the geometric coupling, assuming that the ALE velocity of the fluid domain $\bm u_{\rm ALE}$ and the fluid domain configuration at time $t^{n+1}$ are obtained by extrapolation from previous time steps \citep{nobile2001numerical}. 
The nonlinear convective term of the NS equations is linearized through a first-order semi-implicit Euler scheme, and the structural dynamics is discretized by the implicit BDF1 scheme.

The resulting time-discrete FSI problem is numerically solved by means of the LC-DN-$\alpha$ partitioned algorithm, whose formulation is summarized in Algorithm \ref{alg:DN-alpha-LC}.
At each time step, this loosely-coupled scheme performs one fluid and one structure solution, with suitably modifying the right-hand sides.
Owing to the parameter \al, this scheme features very interesting stability properties even in the presence of significant added-mass.
In Algorithm~\ref{alg:DN-alpha-LC},
$\frac{\delta\bm{v}}{\delta t} = \frac{\partial \bm{v}}{\partial t} + (\bm{u}_{\text{ALE}}\cdot \nabla)\bm{v}$ denotes the ALE time derivative \citep{donea1982arbitrary,quarteroni2017cardiovascular}; 
$\rho_f$ and $\rho_s$ are the fluid and structure densities; $\beta\neq 0$ surrogates the circumferential elastic forces in the 2D case, whereas $\beta=0$ in the 3D case; $\bm{T}_f$ and $\bm{T}_s$ are the fluid and the structure Cauchy stress tensors. 
Notice that, for the sake of notational simplicity, we use the same symbols for the Eulerian and Lagrangian representations of the unknowns.


\begin{algorithm}
\caption{LC-DN-$\alpha$ scheme}
\label{alg:DN-alpha-LC}
{\sl Given suitable initial conditions and boundary conditions on $\partial\Omega_f^t\setminus\Sigma^t$ and $\partial\Omega_s^t\setminus\Sigma^t$, find at $t^{n+1} \, \forall n = 0,1, \dots, T/ \Delta t - 1$, $\bm{u}^{n+1}, p^{n+1}$, ${\boldsymbol{\eta}}^{n+1}$, and $\bm{w}^{n+1}$ by solving in sequence ($n+1$ understood):}
\vspace{1pt}
\begin{enumerate}
    \item Fluid problem equipped with a Dirichlet boundary condition at the interface $\Sigma$: 
    \begin{align}
    \left\{
        \begin{aligned}
            &\rho_f \left(\frac{\bm{u}}{\Delta t} + (\bm{u}^n - \bm{u}_{ALE}^n) \cdot \nabla \bm{u}\right)- \nabla \cdot \boldsymbol{T}_{f} =  \bm{f}_f + \rho_f \frac{\bm{u}^n}{\Delta t} + \notag \\ 
		& \qquad \qquad + (1- \alpha) \left[\rho_f \left(\frac{\bm{u}^{n} -   \bm{u}^{n-1}}{\Delta t} + (\bm{u}^{n} - \bm{u}_{ALE}^{n}) \cdot \nabla \bm{u}^{n}\right) - \nabla \cdot \boldsymbol{T}_{f}^{n}\right]
            &{} \text{in}\,\, \Omega_f,\\[6pt]
            &\nabla \cdot \bm{u} = (1 - \alpha)\nabla \cdot \bm{u}^{n}
            &{} \text{in}\,\, \Omega_f,\\[6pt]
                &\bm{u} = (1 - \alpha)\bm{u}^{n} + \alpha \bm{w}^n
            &{} \text{on}\,\, \Sigma;
        \end{aligned}
    \right.
\end{align}

    \item Structure problem equipped with a Neumann boundary condition at the interface $\Sigma$: 
    \begin{align}
    \left\{
    \begin{aligned}
        &\rho_s \frac{{\boldsymbol{\eta}}}{\Delta t^2} - \nabla \cdot {\boldsymbol{T}}_{s}  + \beta\bm{\eta} = 
        \bm{f}_s + \rho_s\frac{2{\boldsymbol{\eta}}^n - {\boldsymbol{\eta}}^{n-1}}{\Delta t^2} + \\
        & \qquad \qquad \qquad+ (1-\alpha) \left( \rho_s \frac{{\boldsymbol{\eta}}^{n} - 2{\boldsymbol{\eta}}^{n-1} + {\boldsymbol{\eta}}^{n-2} }{\Delta t^2} - \nabla \cdot {\boldsymbol{T}}_{s}^n\right)
        &{} \text{in } {\Omega}_s^0,\notag \\[6pt]
        &\boldsymbol{T}_{s} \cdot \bm{n} =  \boldsymbol{T}_{f} \cdot \bm{n} + (1-\alpha) \left(\boldsymbol{{T}}_{s}^n \cdot \bm{n} - \boldsymbol{T}_{f}^n \cdot \bm{n} \right)
        &{} \text{on}\ \Sigma^0; \notag 
    \end{aligned}
    \right.
    \end{align}
    \item Update $\displaystyle{\bm{w} = \frac{\bm{\eta} - \bm{\eta}^n}{\Delta t}\bigg|_{\Sigma^0}}$.
    \end{enumerate}
\end{algorithm}

\subsection{Stability of the LC-DN-\al scheme}
In this section, we report the main stability property of Algorithm \ref{alg:DN-alpha-LC} obtained in \cite{renzi2026stable}. 
In particular, the following result was obtained for the same benchmark problem described in \cite{causin2005added,gigante2021stability,renzi2026stable}, where the fluid domain consists of a fixed two-dimensional rectangular grid and the structure domain is a one-dimensional membrane with thickness $H_s$.
Moreover, the fluid is modelled as a linear, incompressible, and inviscid fluid, and the structure through a 1D independent rings model \citep{ma1992numerical,quarteroni2000computational}. 
Finally, the structure displacement is assumed to occur only in the direction normal to the interface. 
 

\begin{theorem}\label{th:LC-stability}
The LC-DN-$\alpha$ scheme is unconditionally absolute stable if
\begin{equation}
    0 < \alpha < \overline{\alpha}^{*} = \frac{2\rho_sH_s}{\rho_sH_s + \rho_f\mu_{\text{max}}}.
\end{equation}
\end{theorem}
\noindent where $\mu_{\text{max}}$ is the maximum eigenvalue of the added mass operator \citep{causin2005added, gigante2021stability}.
\newline

The main contribution of introducing the parameter \al is to overcome the instability of the classical DN scheme (recovered with setting $\alpha = 1$) in the presence of large added-mass.
Indeed, when $\rho_sH_s < \rho_f\mu_{\text{max}}$, the critical value defined in Theorem~\ref{th:LC-stability} satisfies $\overline{\alpha}^* < 1$, implying that the standard choice \al$ = 1$ is unstable, consistently with the results of \cite{causin2005added}. 
However, selecting a sufficiently small value of \al restores unconditional absolute stability. 

While this theoretical analysis provides stability bounds for \al in the benchmark problem
and the numerical results presented in \cite{renzi2026stable, renzi2026choice} highlighted the effectiveness of such theoretical 
results even with more complex fluid and structure operators, it does not yield a practical criterion for selecting the optimal value of \al for a given FSI configuration.
Yet, the performance of the LC-DN-\al scheme critically depends on this choice \citep{renzi2026choice}.
Moreover, Theorem \ref{th:LC-stability} shows that \al is not a universal parameter, being dependent on the physical characteristics of the FSI problem (e.g., fluid and structural densities).
This dependence, already evident in the benchmark setting, is expected to become even more pronounced in more complex FSI scenarios, where nonlinear effects, complex geometries, and additional physical parameters may further affect the coupling behavior.

In summary, a proper characterization of the relationship between the problem parameters and the optimal value of \al is essential for the effective application of the LC-DN-\al scheme to general FSI problems. 
However, in its current state, a general model capable of describing this dependence and guiding the selection of \al is lacking.
Here, we aim to overcome this limitation by proposing a data-driven framework for the automatic estimation of \al.

In the following sections, we first introduce the learning framework and then we investigate how the main physical and numerical parameters influence the optimal choice of \al, with the aim of identifying the relevant features governing this relationship. 
The resulting analysis provides the foundation for developing a data-driven model to automatically predict \al in Algorithm \ref{alg:DN-alpha-LC}.

\section{The learning framework}\label{sec:problem_learning}

Let $\bm d$ be the vector of input data related to a single FSI problem instance.
In particular, $\bm d$ collects geometric and configuration data (e.g., domain geometry and structural thickness), physical data (e.g., fluid and solid densities, viscosity, and Lamé parameters), and numerical discretization data (e.g., mesh size and time-step size).
We expect that these data influence the mechanics of the FS coupling, and in turn, the optimal choice of the stabilization parameter \al.

For a given $\bm d$, we consider two formulations of the FSI problem: the monolithic formulation and the one based on the LC-DN-\al scheme (Algorithm \ref{alg:DN-alpha-LC}). 
We denote the monolithic formulation by $P_{\text{M}}(\bm{d})$, which represents the algebraic problem resulting from the discretization of the fully coupled FSI problem:
\[
P_{\text{M}}:\bm d\mapsto\{\bm y^n_{\text{M}}(\bm d)\}_{n=0}^{N_d},
\]
where $\bm y_{\text{M}}$ is the solution vector collecting the discrete unknowns (fluid velocity and pressure, fluid domain and structure displacements) of the monolithic problem at the discrete time instant $n$, with $n = 0, \dots N_d= T/\Delta t$, $T$ the simulation interval, and $\Delta t$ the time-step size.

Equivalently, we denote the LC formulation by $P_{\text{LC}}(\bm{d}, \alpha)$, which represents the algebraic problem resulting from the discretization of the partitioned FSI problem:
\[
P_{\text{LC}}:(\bm d, \alpha) \mapsto \{\bm y_{\text{LC}}^n(\bm d, \alpha)\}_{n=0}^{N_d},
\]
with solution $\bm y_{\text{LC}}$, which depends also on the parameter $\alpha$ of the LC-DN-\al scheme \citep{renzi2026choice}.

For each problem instance defined by $\bm d$, the monolithic solution $\bm y_{\text{M}}(\bm d)$ is taken as the reference.
The objective is to identify the value of \al that ensures the stability of the LC scheme and provides the most accurate approximation of $\bm y_{\text{M}}(\bm d)$.
Accordingly, we define the optimal stabilization parameter as
\begin{equation}\label{eq:optimal_alpha}
    \alpha^*(\bm{d}) = \operatorname*{arg\,min}_{\alpha\in\mathcal{I}} \left \{\max_n \left\|\bm{y}_{\text{M}}^n(\bm{d}) - \bm{y}_{\text{LC}}^n(\bm{d},\alpha)\right\|^2\right\},
\end{equation}
where $\mathcal{I}$ denotes the set of stable values of \al for the considered problem instance, and $\|\cdot\|$ is a suitable norm, such as the $H^1$ norm.

Regardless of the optimization method used to solve \eqref{eq:optimal_alpha}, computing $\alpha^*(\bm{d})$ requires multiple solutions of the LC-DN-\al scheme for a given $\bm{d}$.
Thus, we aim to approximate the mapping from the problem-instance parameters to the optimal stabilization parameter,
\[
\mathcal A(\bm d) :\bm d\mapsto\alpha^*(\bm d),
\]
using machine-learning methods. Specifically, we introduce a feed-forward neural network
\[
f_{\bm{\theta}}:\bm{d}\mapsto\widehat{\alpha},
\]
where $\bm{\theta}$ denotes the trainable parameters of the network (weights and biases) and $\widehat{\alpha}$ is the estimated stabilization parameter. 

The training dataset is generated in an offline phase, where for each problem instance $\bm{d}_i$ the monolithic formulation is first solved to obtain the reference solution $\{\bm{y}_{{\text M},i}^n\}_{n=0}^{N_d}$.
Then, the LC-DN-\al problem is solved for a set of candidate $\alpha_{ij}$, with each value yielding a corresponding $\{\bm{y}_{{\text LC},i}^n(\bm{d}_i,\alpha_{ij})\}_{n=0}^{N_d}$.
Then, for each $\bm{d}_i$, the subset $\mathcal{I}_i$ of stable $\alpha_{ij}$ is retained, and the optimal value associated with $\bm d_i$ is computed as in \eqref{eq:optimal_alpha}.
This offline procedure leads to the supervised learning dataset
$\mathcal{D}=\{(\bm{d}_i,\alpha_i^*)\}_{i=1}^{N}$, where $N$ is the number of instances.

In general, we expect that not all the geometric, physical, and numerical data included in $\bm{d}$ have the same influence on $\alpha^*$. 
Additionally, since a neural network's performance strongly depends on the choice of input features used to represent the problem instance, it is important to identify an effective input representation for the learning model, avoiding redundancies. 

For these reasons, in Section \ref{sec:numerical_results1}, we present an analysis of the dependence of $\alpha^*$ on the different components of $\bm d$ before defining the neural network architecture. 
Specifically, we first introduce the baseline FSI problem (Section \ref{subsec:baseline_pb}) upon which we generate a dataset of problem instances by varying some physical and numerical data, that is the structure's Young modulus $E$ and density $\rho_s$, the mesh size $h$, and $\Delta t$, and we compute the associated $\alpha^*$ (Section \ref{subsec:dataset}).
We then analyze how $\alpha^*$ depends not only on the data used to generate the problem instances dataset but also on selected descriptors that characterize the algebraic structure of the resulting discrete FSI problem (Section \ref{subsec:features}).
This analysis allows us to identify three candidate input representations for the neural network ${\bm{d}}_I$, ${\bm{d}}_{II}$, and ${\bm{d}}_{III}$ defined as subsets of the components of the original data vector $\bm d$.
Finally, we train a separate neural network $f_{\bm{\theta}, k}$ for each representation ${\bm{d}}_k$, with $k = I, II, III$, and we compare their performance to select the best-performing model (Section \ref{subsec:NNtraining}), which is used for all subsequent analyses in the remainder of the paper.

\section{Results: problem setup and neural network architecture}\label{sec:numerical_results1}
In this section we introduce the problem setup and the neural network architecture.
\subsection{The baseline problem}\label{subsec:baseline_pb}
\begin{figure}
  \centering
  \includegraphics[width=0.6\textwidth]{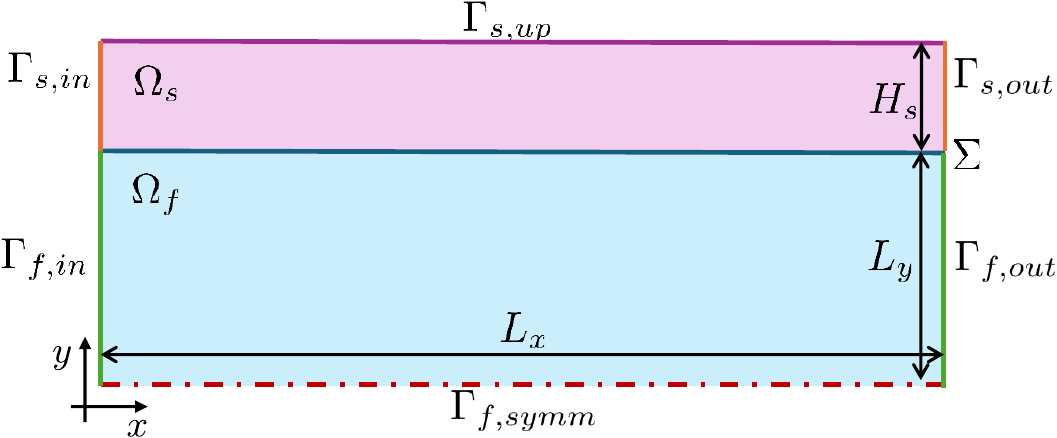}
  \caption{Fluid and structure domains for the baseline problem.}
  \label{fig:domain}
\end{figure}

For the neural network training, we consider the 2D FSI problem defined on a rectangular domain, representing a longitudinal section of a cylindrical blood vessel, as illustrated in Fig. \ref{fig:domain}.
The fluid is modeled as an incompressible Stokes flow, while the vessel wall is described as a linear elastic structure. The geometrical, fluid, and structural parameters defining the baseline problem are summarized in Table \ref{tab:baselinepb}.

\begin{table}[ht]
\centering
\caption{Geometrical, fluid, and structural parameters of the baseline FSI problem.}
\label{tab:baselinepb}
\begin{tabular}{lll}
\hline
\textbf{Parameter} & \textbf{Symbol} & \textbf{Value} \\
\hline
Fluid domain length & $L_x$ & $6\,cm$ \\
Fluid domain height & $L_y$ & $1\,cm$ \\
Structure thickness & $H_s$ & $0.1\,cm$ \\
Fluid dynamic viscosity & $\mu$ & $0.035 \,Poise$ \\  
Fluid density & $\rho_f$ & $1.0\, g/cm^3$ \\
Structure density & $\rho_s$ & $1.1 \,g/cm^3$    \\
Young's modulus & $E$ & $1.3\times10^6\, dyne/cm^2$ \\
Poisson's ratio & $\nu$ & $0.3$ \\
\hline
\end{tabular}
\end{table}

Regarding boundary conditions, for the Stokes problem, we impose a pressure profile at the inlet $\Gamma_{f,in}$:
\begin{equation}\label{eq:pressure_in}
    p_{in} = 
    \begin{cases}
        4\times 10^4\left(1 - \cos\left( \frac{\pi t}{0.05}\right)\right) \, dyne/cm^2 \qquad t\le 0.1s, \\
        0 \, dyne/cm^2 \qquad\qquad\qquad\qquad\qquad\qquad 0.1s<t<T,
    \end{cases}
\end{equation}
homogeneous Neumann conditions, $\bm{T}_f\cdot \bm{n}= \bm{0}$
\footnote{We consider here a reflecting boundary condition. An absorbing boundary condition could also be considered; however, since the right-hand side strongly depends on both the data and time instants, the corresponding training procedure proved to be more difficult than in the reflecting case. 
We nevertheless assess the effectiveness of the network in the absorbing case as well; see Section \ref{sec:numerical_results}.}, at the4 outlet $\Gamma_{f,out}$, and symmetry conditions, $(\bm{T}_{f}\cdot \bm{n})\cdot \bm{t} = 0 \text{ and } \bm{u}\cdot \bm{n}= 0$, on $\Gamma_{f,symm}$.
For the structure problem, we imposed null normal displacement, $({\bm{T}}_s \cdot \bm{n})\cdot {\bm t} = 0 \text{ and } {\bm{\eta}}\cdot \bm{n}= 0$, on $\Gamma_{s,in}^0 \cup \Gamma_{s,out}^0$, and free stresses $\bm{T}_s\cdot\bm{n} = \bm{0}$ on $\Gamma_{s,up}^0$.
Regarding the numerical settings, we adopt the time discretization schemes detailed in Section \ref{subsec:formulation}, with time step $\Delta t$, and $\mathbb{P}_1-iso\mathbb{P}_2/\mathbb{P}_1$ Finite Elements in space for the fluid and $\mathbb{P}_1$ for structure variables, with characteristic mesh size $h$.

\subsection{Generation of the learning dataset}\label{subsec:dataset}
Starting from the baseline problem, a dataset of $N$ FSI configurations, with $N=361$, is generated by applying controlled variations to key data that affect the stabilization behavior. 
In particular, we vary $h,\Delta t,E, \text{ and }\rho_s$ within prescribed ranges, which are reported in Table \ref{tab:ranges}.

\begin{table}[htpb]
\caption{Ranges of variability of the data defining the FSI configurations employed to generate the learning dataset.}\label{tab:ranges}
\begin{center}
\renewcommand{\arraystretch}{1.5}

\begin{tabular}{l|l}
\toprule
 Data & Range \\ 
\midrule
 $h$   & $(0.0625, 0.500)cm$\\
 $\Delta t$ & $(3.125\times 10^{-5}, 1.000 \times 10^{-3})s$\\
 $E$ & $(1.000 \times 10^{5}, 1.000\times 10^{7})dyne/cm^2$\\
 $\rho_s$ & $(0.500, 1.500)g/cm^3$\\
\bottomrule
\end{tabular}
\end{center}
\end{table}

\noindent To ensure an approximately uniform distribution of the samples within the four-dimensional parameters' space, we adopted an optimized Latin Hypercube Sampling strategy based on the minimization of the centered discrepancy \citep{mckay2000comparison, damblin2013numerical}. 

For each sampled configuration, the corresponding optimal stabilization parameter $\alpha^*$ is computed following the procedure described in Section \ref{sec:problem_learning}.
The simulations are carried out using a 2D Finite Element Matlab code developed at LaBS, Dipartimento di Chimica, Materiali e Ingegneria Chimica of Politecnico di Milano.

\subsection{Input space characterization}\label{subsec:features}

In this section, we explore how the optimal value of the stabilization parameter $\alpha^*$ depends on the data of the FSI problems, to select the most representative set of input features for the learning model.  

The candidates include both physical (i.e. $\rho_s$ and $E$) and numerical (i.e. $h$ and $\Delta t$) quantities.
In addition, we consider some descriptors of the algebraic structure of the FSI matrix system, such as the Frobenius norm of the fluid and structure matrices $\bm{A}_f \text{ and } \bm{A}_s$ and their inverses, and some algebraic descriptors of the FS interface coupling, such as the Frobenius norm of the Schur complement $\bm{\Sigma}$ associated with the coupled FSI system and the Frobenius norm of fluid and structure contributions of the Schur complement, denoted as $\bm{\Sigma}_f$ and $\bm{\Sigma}_s$, respectively. 
We define {\sl candidate input data} the set $\{\rho_s, E, h, \Delta t, \|\bm{A}_f\|_F, \|\bm{A}_f^{-1}\|_F, \|\bm{A}_s\|_F, \|\bm{A}_s^{-1}\|_F, \|\bm{\Sigma}\|_F, \|\bm{\Sigma}_f\|_F, \|\bm{\Sigma}_s\|_F\}$.

\begin{figure}
  \centering
  \includegraphics[width=1\textwidth]{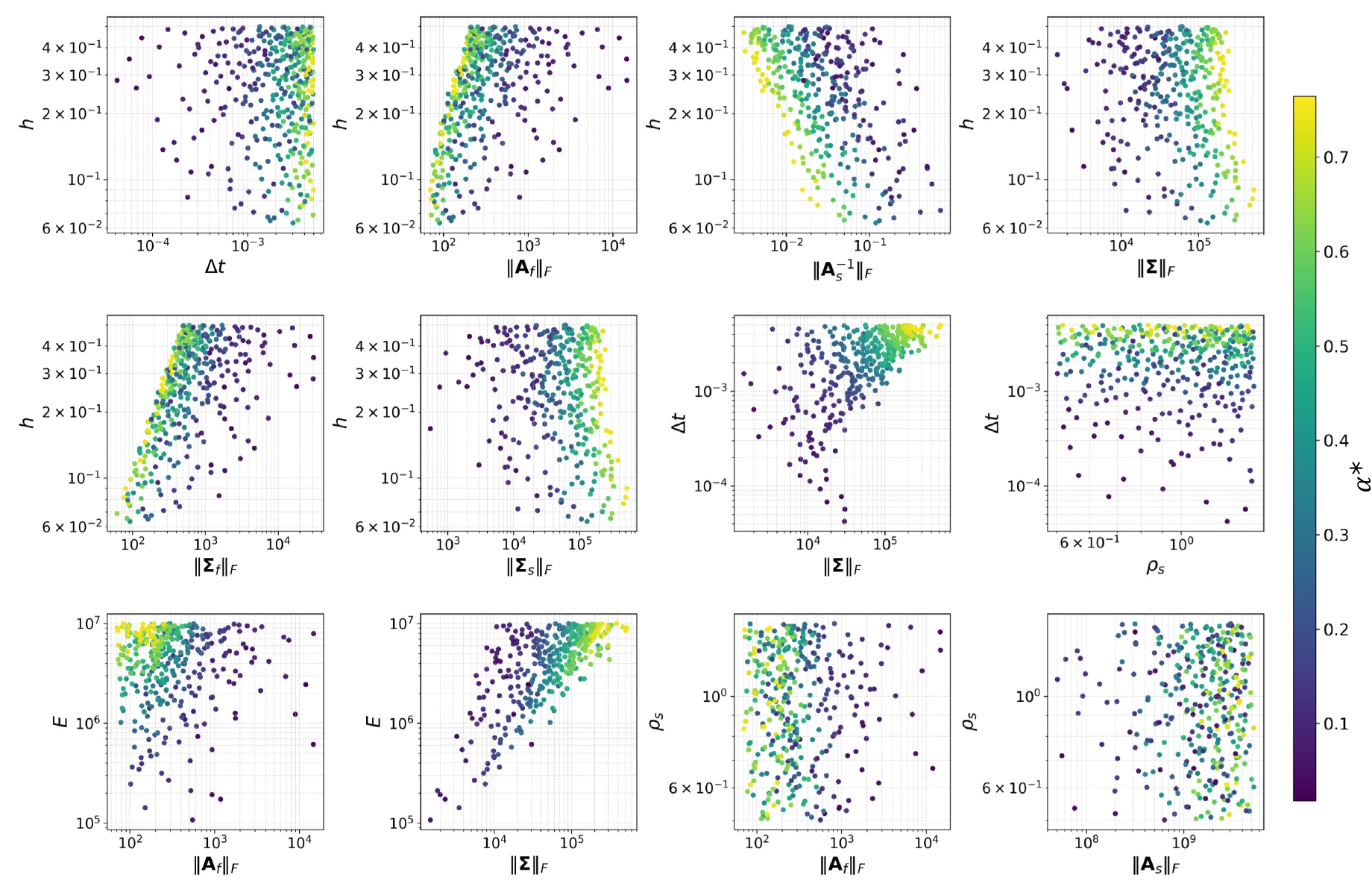}
  \caption{Scatter plots representing the pairwise distributions of candidate input data of the FSI problems.
  Each point corresponds to a specific FSI problem, and it is colored according to the corresponding optimal value of $\alpha$.}
  \label{fig:scatter}
\end{figure}

To qualitatively assess the relationship between the candidate input data and the optimal stabilization parameter \al, we analyze their pairwise distributions through scatter plots.
Fig. \ref{fig:scatter} reports a selection of the most representative data pairs, with each sample colored according to the corresponding value of $\alpha^*$. 
Through this analysis, we identify several informative pairs, which are characterized by clustered distributions and nonlinear patterns in the values of $\alpha^*$, e.g., the pairs $(E,\|\bm{\Sigma}\|_F)$ and $(\Delta t,\|\bm{\Sigma}\|_F)$. 
We observe that the relationship between the data and 
$\alpha^*$ is generally nonlinear, thus motivating the use of a nonlinear regression model. 

\begin{figure}
  \centering
  \includegraphics[width=0.95\textwidth]{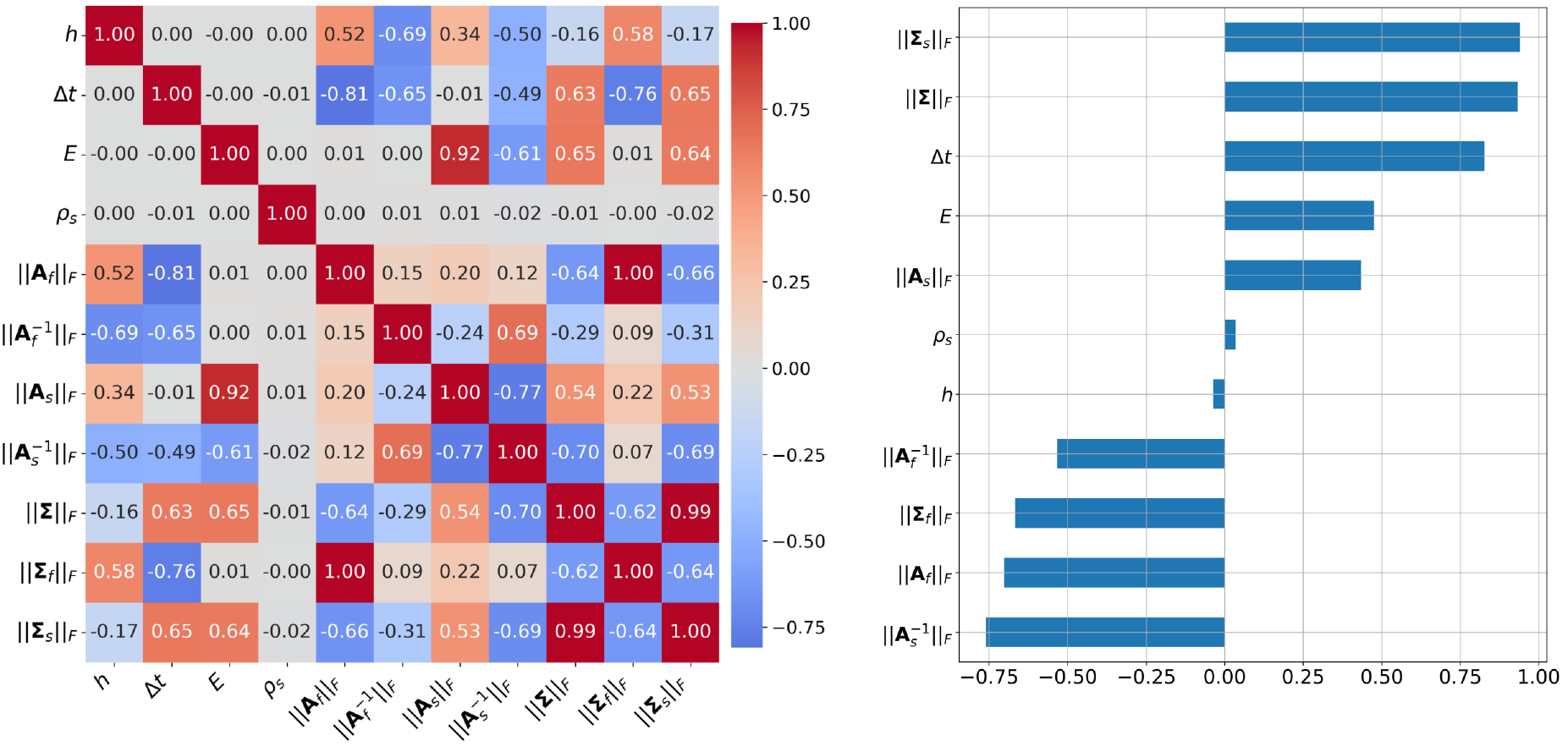}
  \caption{Left: Spearman correlation among the candidate input data. Right: Spearman correlation between the candidate input data and $\alpha^*$.}
  \label{fig:correlations}
\end{figure}

To further investigate the relevance and redundancy of the candidate input data, we compute the pairwise correlations among the candidate input data (Fig. \ref{fig:correlations}-left) and between each data and $\alpha^*$ (Fig. \ref{fig:correlations}-right).
The results reveal that the algebraic descriptors exhibit strong positive or negative correlations with $\alpha^*$. However, these quantities are also highly correlated one another and with $E$ and $\Delta t$.
Moreover, $E$ and $\Delta t$ show a strong positive correlation with $\alpha^*$.
Surprisingly, $\rho_s$ and $h$ display only very weak correlations with $\alpha^*$ and with any other feature. 

Based on these arguments, we identify three candidate input sets: 
\begin{itemize} 
    \item ${\mathcal{D}}_{I} = \{{\bm d}_{I,i}\}_{i=1}^N$, with ${\bm d}_{i,I} = \{\Delta t_i, E_i, \rho_{s,i}\}$. This input set includes the physical and numerical parameters that most strongly influence the value of $ \alpha^*$.
    We also added $\rho_s$ as it is intrinsically related to the added-mass mechanism and is therefore expected to be linked to $\alpha^*$, although its influence is less evident in this simplified problem.    

    \item ${\mathcal{D}}_{II} = \{{\bm d}_{ II,i}\}_{i=1}^N$, with ${\bm d}_{i, II} =  \{\Delta t_i, \|\mathbf{A}_s\|_{F,i}, \|\mathbf{A}_f\|_{F,i}\}$. In this input set, we include the Frobenius norms of the two matrices that exhibit strong, opposite correlations with $\alpha^*$. $\Delta t$ is also retained since the scatter plots indicate an almost quadratic relationship between $\Delta t$ and $\alpha$.

    \item ${\mathcal{D}}_{III} = \{{\bm d}_{III, i}\}_{i=1}^N$, with ${\bm d}_{i, III} = \{\Delta t_i, E_i, \rho_{s,i}, \|\mathbf{\Sigma}\|_{F,i}\}$. This third input set extends ${\bm d}_{I}$ by including $\|\mathbf{\Sigma}\|_F$, with the aim of providing the model with additional information about the coupling effects governing the added-mass behavior.
    
\end{itemize}

Notice that, for construction, the values of $\Delta t_i $ does not depend on $k= I,II, III$, and same holds true for all the other input parameters.

\subsection{Selection and training of the neural network model}\label{subsec:NNtraining}
Three feed-forward neural networks (FFNNs) were investigated to predict the stabilization parameter $\alpha^*$.
The three models, denoted as $f_{\theta,I}$, $f_{\theta, II}$, and $f_{\theta,III}$, share the same architecture and training procedure and differ only in the set of input data adopted to characterize the FSI problem, namely ${\mathcal{D}}_I, {\mathcal{D}}_{II}, \text{ and }{\mathcal{D}}_{III}$, respectively.

For each model, the learning dataset ${\mathcal{D}}_k$, with $k= I, II, III$, consists of $N=361$ FSI configurations ${\bm{d}}_{k,i}$, each associated with the corresponding optimal stabilization parameter $\alpha^*_i$. 
The dataset is randomly divided into a {\sl training set} and a {\sl test set}, containing 70\% and 30\% of the samples, respectively. 
The training set is used to optimize the network parameters (weights and biases), whereas the test set is used to assess the predictive capability of the trained models on unseen configurations.

All networks consist of a fully connected feed-forward architecture with a single hidden layer composed of 10 neurons.
The hidden layer employs the Rectified Linear Unit (ReLU) activation function \citep{sharma2017activation}, while the output layer is linear to provide a continuous estimate of $\alpha^*$.
Prior to training, all input features are normalized in order to improve the optimization process and ensure comparable input scales. The network parameters are optimized using the Adam optimizer \citep{kingma2014adam} with a learning rate of $10^{-3}$. 
As an objective function to minimize during the FFNN training we adopted the Huber loss function
with parameter $\delta=0.5$. 
This function combines the sensitivity of the mean squared error to small prediction errors with the robustness of the mean absolute error to larger prediction errors \citep{huber1992robust}. 

\begin{figure}
  \centering
  \includegraphics[width=1\textwidth]{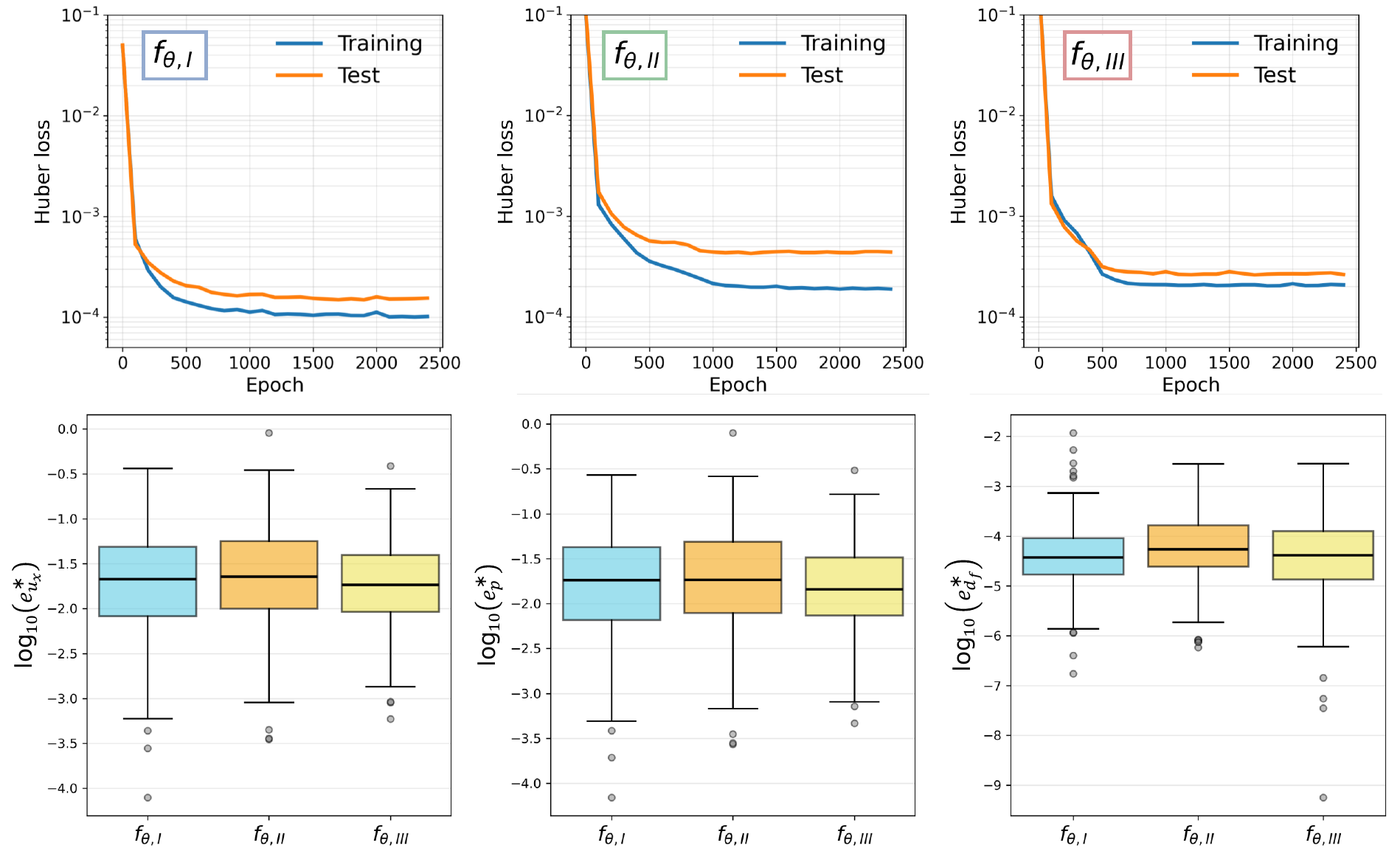}
  \caption{Top row: convergence histories of the training and test losses for the three models $f_{\bm {\theta},I},\, f_{\bm {\theta},II}, \, \text{and}\, f_{\bm {\theta},III}$. Bottom row: boxplots of the logarithm of the relative differences in the test set, computed as \eqref{eq:error_measure}, obtained with the three models. 
  These differences are evaluated for the $x$ component of the fluid velocity $u_x$, the fluid pressure $p$, and the fluid domain diameter $d_y$ (from left to right).}
  \label{fig:NNchoice}
\end{figure}

We train each FFNN for 2500 epochs using mini-batch gradient descent with a batch size of 32.
The convergence histories, in terms of the Huber loss, of the training and test losses are reported in Fig. \ref{fig:NNchoice}-top for the three models. In all cases, both losses decrease steadily and reach a stable plateau, indicating successful convergence, and remain close to each other, indicating good generalization.

In addition, we also compare the performance of the three FFNN by considering, for the test set, the accuracy of the resulting LC FSI solution when the predicted value of $\widehat{\alpha}$ is used in the numerical solver, with respect to the LC FSI solution obtained using the corresponding $\alpha^*$.
In particular, we define, for a generic scalar discrete variable $v$ an error measure as: 
\begin{equation}\label{eq:error_measure}
    e^*_{v} = \frac{\max_{0<n<N_t} \frac{h^2}{|\Omega^0_{h,\bullet}|}\sum_{j}\left(\overline{v^{*,n}}_{j} - \overline{\widehat{v^n}}_{j}\right)^2}{\max_{0<n<N_t} \frac{h^2}{|\Omega^0_{h,\bullet}|}\sum_{j}\left(\overline{v^{*,n}}_{j}\right)^2},
\end{equation}
where, $N_t$ is the total number of time steps, $|\Omega_{h,\bullet}| \text{ with } \bullet=f,s$ is a measure of the fluid/structure domain, $j$ is the index of the degrees of freedom along the $x$ coordinate.
With ${v^*}$ we indicate the LC-DN-\al solution obtained using $\alpha^*$, and $\widehat{v}$ is the LC-DN-\al  solution obtained using $\widehat{\alpha}$. 
The overline denotes averaging along $y$ coordinate.
Moreover, in view of the results presented in Section \ref{sec:numerical_results}, we introduce here the error with respect to the monolithic solution: 
\begin{equation}\label{eq:error_measure_mono}
    e^{M}_{v} = \frac{\max_{0<n<N_t} \frac{h^2}{|\Omega^0_{h,\bullet}|}\sum_{j}\left(\overline{v^{M,n}}_{j} - \overline{\widehat{v^n}}_{j}\right)^2}{\max_{0<n<N_t} \frac{h^2}{|\Omega^0_{h,\bullet}|}\sum_{j}\left(\overline{v^{M,n}}_{j}\right)^2},
\end{equation}
where $v^{M}$ represents the monolithic solution.

Fig. \ref{fig:NNchoice}-bottom reports with boxplots the errors $e^*_v$ obtained, for the test set, for the $x$ component of the fluid velocity $u_x$, the fluid pressure $p$, and the fluid domain diameter $d_f$.
The large overlap of the interquartile ranges indicates that the error distributions are of the same order of magnitude for the vast majority of the tested FSI configurations.
In summary, although $f_{\theta,I}$ consistently achieves the lowest training/test loss and slightly better error distributions (lowest median), all three models produce stabilization parameters that lead to comparable solution accuracy.
These results suggest that, for this baseline problem, the choice of input features has a limited influence on the final accuracy of the LC solution. 
For this reason, in the remainder of the paper, we adopt $f_{\theta,I}$ for all the subsequent numerical experiments as this choice does not involve further computations.

\section{Results: robustness of the neural network}\label{sec:numerical_results}

In this Section, we present the results of the numerical experiments performed using the value of the stabilization parameter estimated by the FFNN $f_{\bm{\theta}, I}$, with the aim of assessing its robustness with respect i) to the test set and ii) to changes in the baseline problem.
In what follows, we denote as LC-DN-${\alpha}^*$ and LC-DN-$\widehat{\alpha}$ the LC-DN-\al schemes obtained using $\alpha=\alpha^*$ and $\alpha=\widehat{\alpha}$, respectively.

\subsection{Robustness with respect to the test set}\label{subsec:consistency}
\begin{figure}
  \centering
  \includegraphics[width=0.98\textwidth]{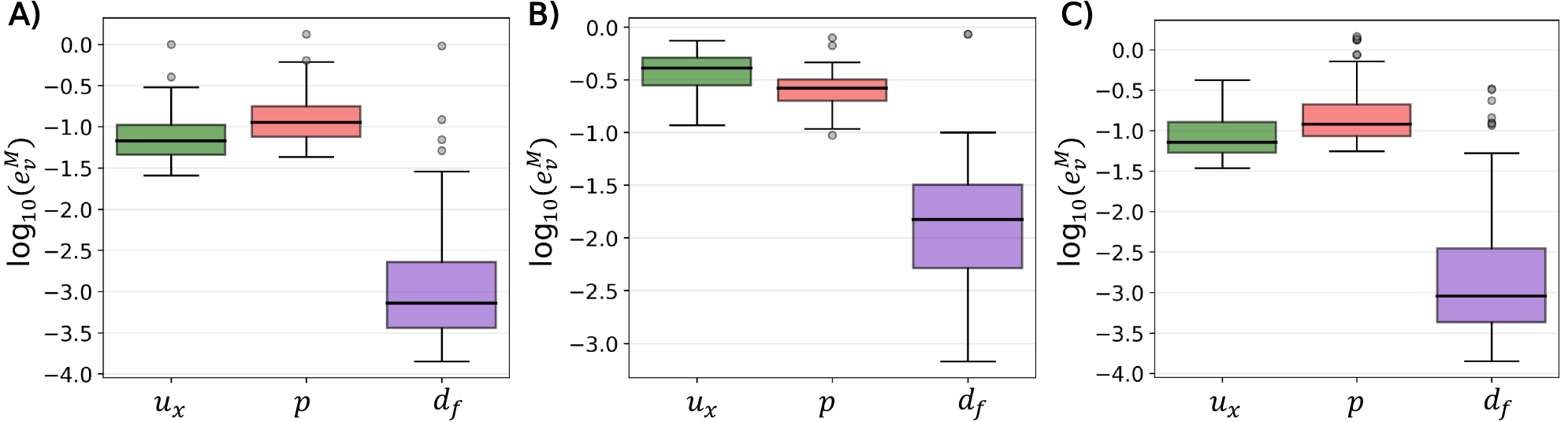}
  \caption{Boxplots of the logarithm of the relative difference between the LC-DN-$\widehat{\alpha}$ solution and the monolithic solution computed as \eqref{eq:error_measure_mono} for the $x$ component of the fluid velocity $u_x$, the fluid pressure $p$, and the fluid domain diameter $d_y$. A) Robustness with respect to the test set; B) Robustness with respect to changes in the BP: P-outflow; C) Robustness with respect to changes in the BP: P-stenosis.}
  \label{fig:errors_network}
\end{figure}

To assess the consistency of the predicted stabilization parameter $\widehat{\alpha}$, we compare the LC solution obtained with $\widehat{\alpha}$ with the solution computed with a monolithic solver.
The comparison is performed across all scenarios in the test set, i.e. that were not used during the neural network's training phase (see Section \ref{subsec:NNtraining}).
The monolithic solution is taken as the reference, and the difference between the two solutions is evaluated in terms of the error measure \eqref{eq:error_measure_mono}.
The results are reported in Fig. \ref{fig:errors_network}-A.

Comparing Fig. \ref{fig:errors_network}-A with Fig. \ref{fig:NNchoice}-bottom row, we observe that the upper whiskers of the discrepancies between LC-DN-$\widehat{\alpha}$ solution and the monolithic solution are comparable to those of the prediction errors associated with $f_{\bm{\theta}, I}$ (light-blue boxplots).
This suggests that the largest discrepancies with respect to the monolithic solution are not significantly affected by the network prediction, while the systematic shift in the median, which is about one order of magnitude larger, is mainly attributable to the intrinsic difference between the loosely-coupled formulation LC-DN-\al and the monolithic solver.  

\begin{figure}
  \centering
  \includegraphics[width=0.98\textwidth]{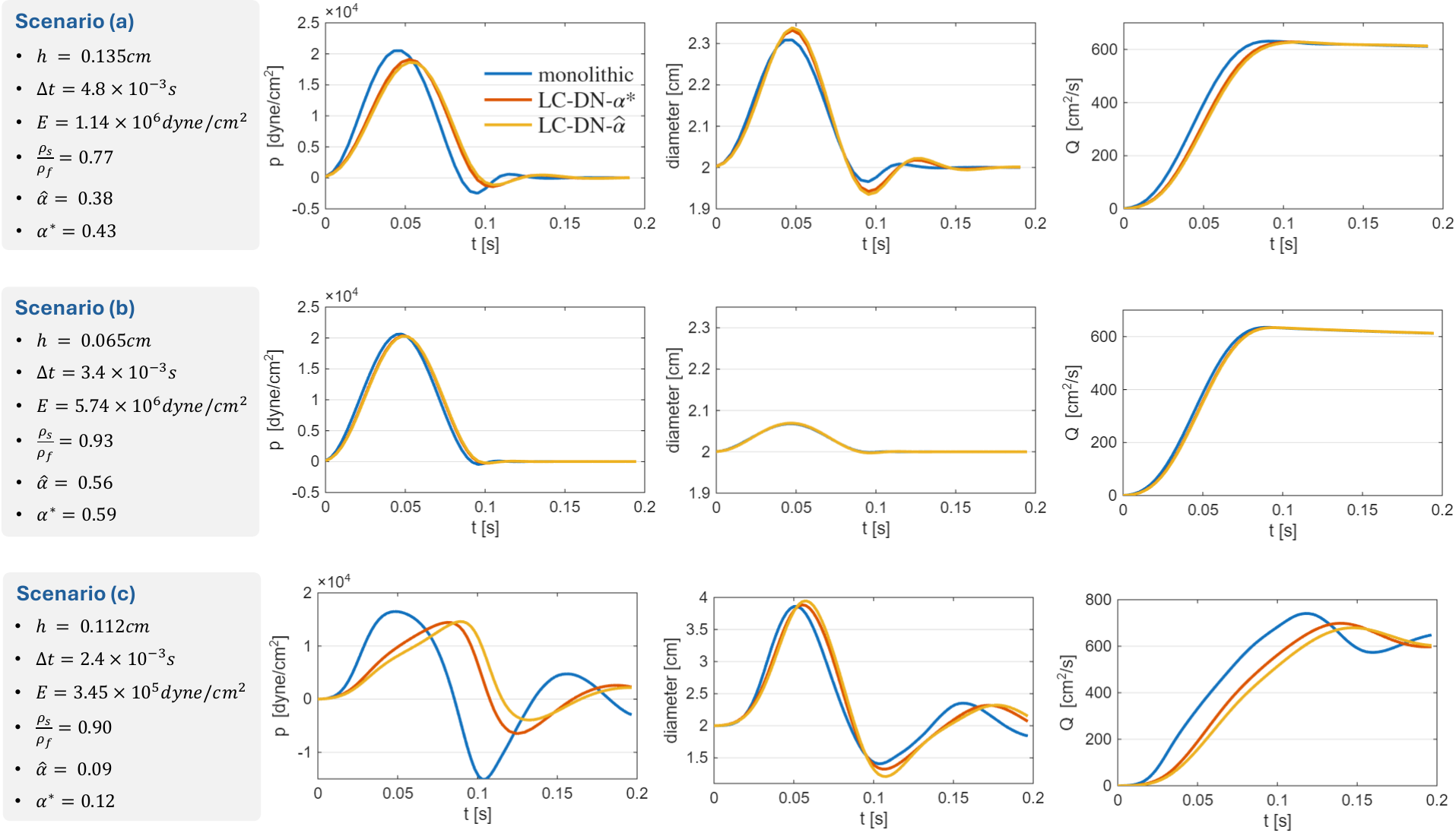}
  \caption{Comparison between the evolution over time of the monolithic, LC-DN-${\alpha}^*$, and LC-DN-$\widehat{\alpha}$ solutions (blue, orange and yellow lines, respectively) extracted at the mid-section of the domain ($x=3cm$) for three representative scenarios.
  The gray panels on the left show the specifics of the scenarios. From the left: the first column shows the fluid pressure averaged along the y-direction, the second the diameter of the fluid domain, and the third the fluid flow. Robustness with respect to the test set.}
  \label{fig:curves_cons}
\end{figure}

\begin{figure}
  \centering
  \includegraphics[width=0.98\textwidth]{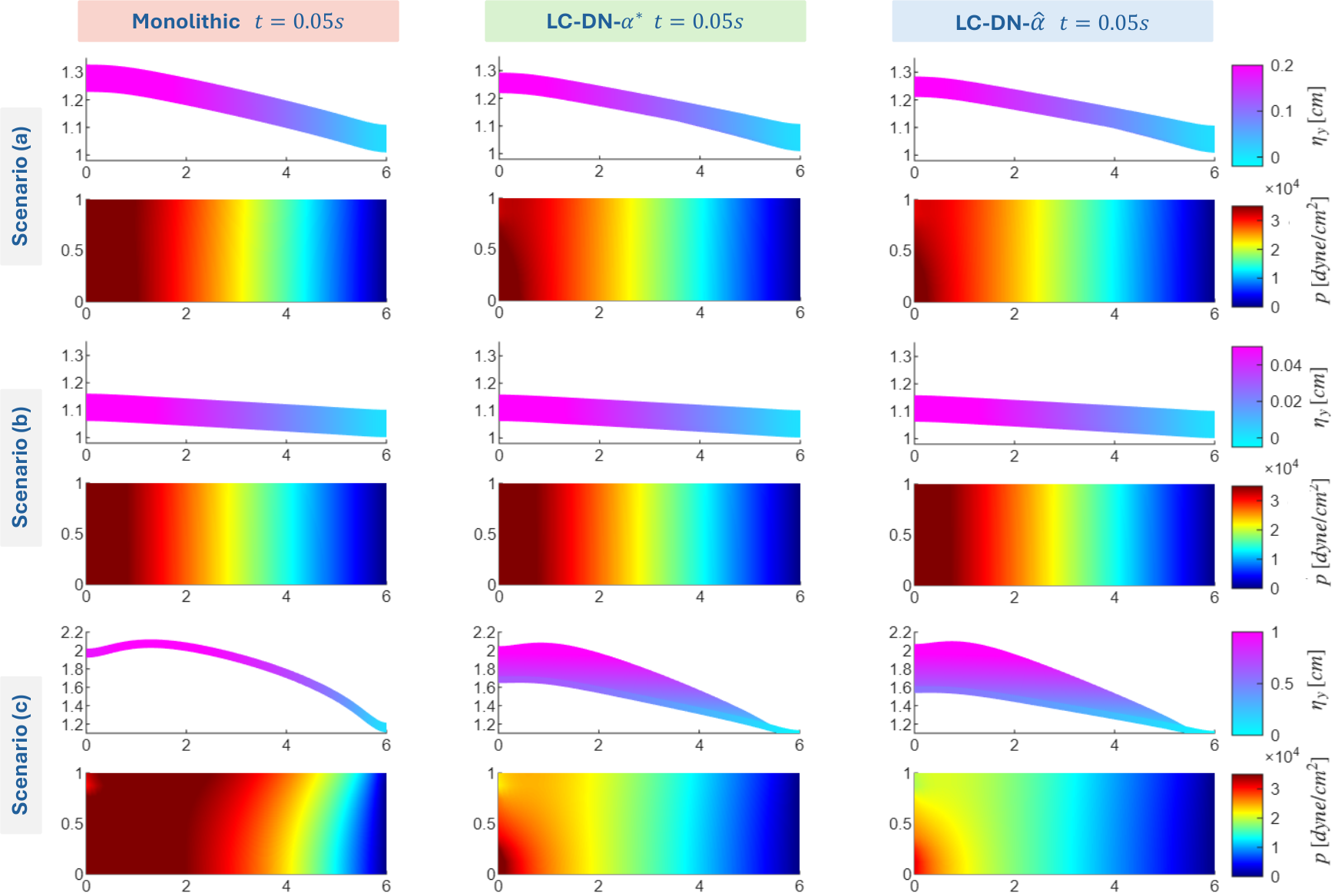}
  \caption{Plots of the LC-DN-$\widehat{\alpha}$ and the monolithic solutions extracted at two time instants, for the same three FSI scenarios detailed in Fig. \ref{fig:curves_cons}. 
  For each scenario, the structure displacement along x and the fluid pressure (top and bottom rows, respectively) are reported. Robustness with respect to the test set.}
  \label{fig:plots_cons}
\end{figure}

This is confirmed by the results of Fig. \ref{fig:curves_cons}, where we report a comparison of the averaged solutions over time obtained with the three solvers (namely the monolithic, LC-DN-$\alpha^*$, and LC-DN-$\widehat{\alpha}$) for three different FSI scenarios randomly selected from the test set. 
We observe that LC-DN-$\widehat{\alpha}$ solution well approximates in any case LC-DN-$\alpha^*$ one, also in scenario (c), where the explicit solutions exhibit larger discrepancies with respect to the monolithic one.
This is also reflected in Fig. \ref{fig:plots_cons}, showing the spatial distribution of the structure displacement and fluid pressure across all three scenarios.
Better accuracy is achieved in scenarios (a) and (b) characterized by higher structure stiffness than scenario (c).
In general, the LC solution obtained with $\widehat{\alpha}$ is close to the one obtained with the optimal value of $\alpha$, for all three scenarios, thus highlighting the good predictive capabilities of the network. 
 
In summary, since the test set consists of unseen configurations, these findings prove the generalization capability of the proposed model.

\subsection{Robustness with respect to changes in the baseline problem}\label{subsec:robustness}
In this section, we consider the baseline problem BP reported in Section \ref{subsec:baseline_pb} and some modified problems, namely {\sl P-inflow}, obtained by changing the inlet condition, {\sl P-outflow}, obtained by changing the outlet condition, and {\sl P-stenosis}, obtained by introducing a stenosis in the original geometry. 
The idea, here, is to assess the accuracy of the LC-FSI solution obtained with $\widehat{\alpha}$, given by the network trained in BP (as described in Section \ref{subsec:NNtraining}), when applied to problems P-inflow, P-outflow, and P-stenosis, i.e., without any ad-hoc training of the network.
The resulting solutions are compared with those obtained using the reference monolithic solver.

\begin{figure}
  \centering
  \includegraphics[width=0.75\textwidth]{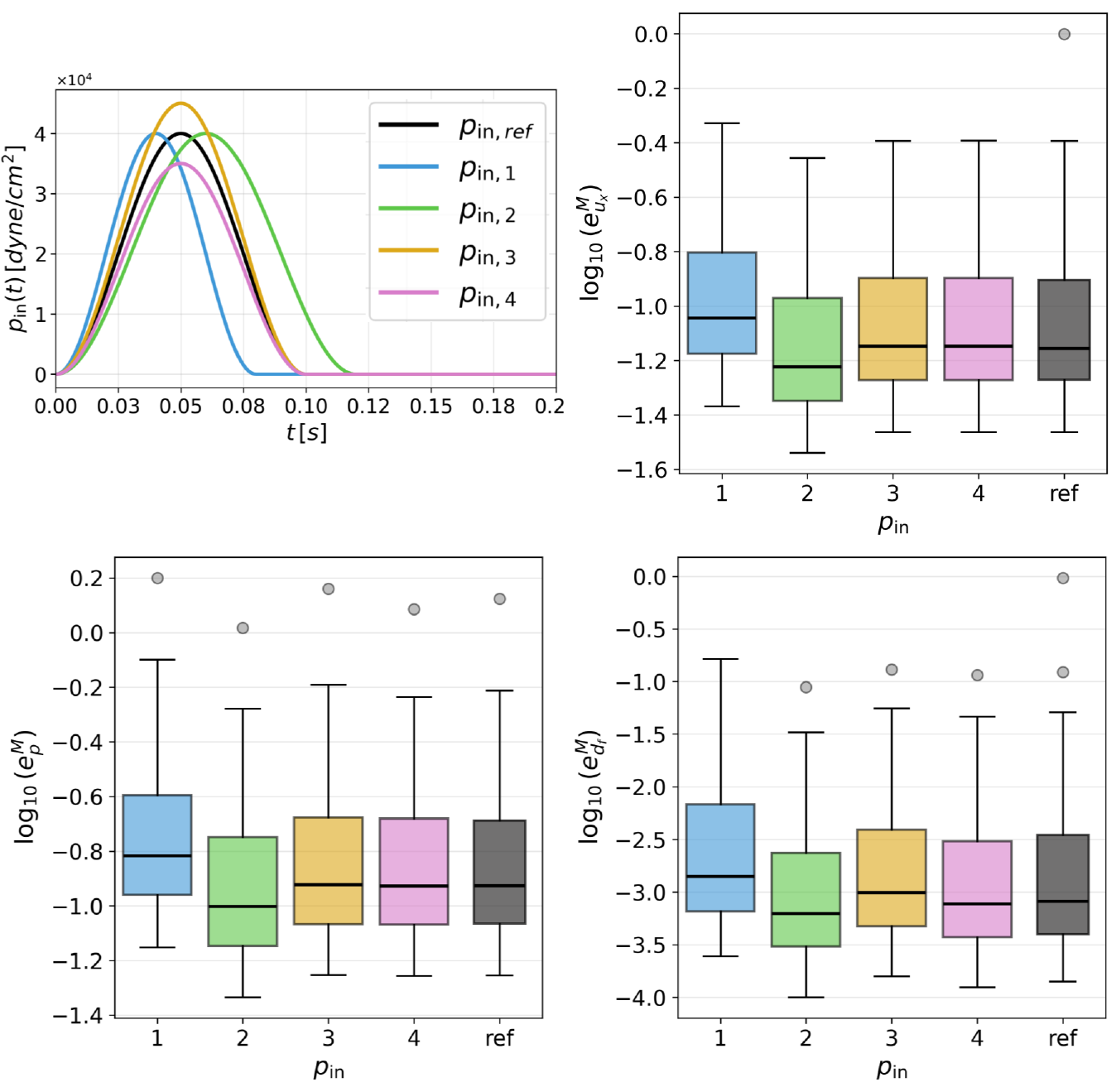}
  \caption{Top left figure: plot of the considered pressure waveform applied as boundary condition at fluid inlet. 
  In the remainder: boxplots of the logarithm of the difference between the LC-DN-$\widehat{\alpha}$ solution and the monolithic solution, computed as \eqref{eq:error_measure_mono}, obtained for all the inlet condition variations. Robustness with respect to changes in the BP: P-inflow.}
  \label{fig:fxin}
\end{figure}

\paragraph{\textbf{P-inflow: robustness to inlet condition}}
In this test, we introduce small variations to the inlet pressure profile \eqref{eq:pressure_in} in terms of period and amplitude of the cosine: the period is modified by approximately $20\%$ (see $p_{\text{in},1}$ and $p_{\text{in},2}$ in the upper-right plot of Fig. \ref{fig:fxin}); the amplitude is varied by approximately $6\%$ (see $p_{\text{in},3}$ and $p_{\text{in},4}$ in the same plot).

Almost all the simulations lead to stable solutions ($100\%$ of the cases with $p_{\text{in},1}$, $p_{\text{in},2}$, and $p_{\text{in},4}$  and $99.8\%$ with $p_{\text{in},3}$).
In Fig. \ref{fig:fxin}, we report the differences between these LC-DN-$\widehat{\alpha}$ solutions and the monolithic reference.
Although $p_{\text{in}, 1}$ is characterized by slightly worse performance, the error distributions are comparable across all modified inlet conditions, with median error values consistently below $10\%$ for $u_x$, $16\%$ for $p$, and below $0.3\%$ for $d_f$.
Hence, we do not observe a significant deterioration in the accuracy with respect to the baseline simulations for which $\widehat{\alpha}$ was estimated (black boxplot in Fig. \ref{fig:fxin}). 

Overall, these results indicate that moderate perturbations of the inlet pressure profile, both in amplitude and period, do not significantly affect the accuracy of the LC-DN-$\widehat{\alpha}$ solutions, suggesting that the network prediction for BP remains effective even when the inlet flow conditions vary slightly.

\paragraph{\textbf{P-outflow: robustness to outlet condition}}
In this test, we replace the homogeneous Neumann condition originally prescribed at the fluid outlet for the training of the network with an absorbing boundary condition \citep{nobile2008effective,nobile2014inexact} in order to reduce non-physical pressure-wave reflections caused by the artificial truncation of the computational domain.

We find that only $44.4\%$ of the FSI simulations remain stable, and the accuracy of the stable solutions is appreciably worse than the baseline results as highligthed by Fig. \ref{fig:errors_network}-B and Fig. \ref{fig:errors_network}-A, respectively. 
Specifically, the median of the relative errors increases up almost to $50\%$ for $u_x$, $30\%$ for $p$, and up to $3\%$ for $d_f$. 
Notice that this behaviour is due to the use of a network trained for the reflective case in a scenario with absorbing conditions.

The reduced effectiveness of the predicted stabilization parameter in this test can be attributed to the modification of the fluid dynamic response induced by the new outlet condition. 
We conclude that, unlike the previous test, altering pressure-wave propagation modifies the fluid-structure interaction pattern, thus affecting the stability and accuracy properties of the partitioned scheme.
Thus, new training is recommended for the absorbing case.


\paragraph{\textbf{P-stenosis: robustness to geometry}}
In this test, we modify the computational geometry by introducing a localized narrowing in the central section of the fluid domain, thereby mimicking the configuration of a stenotic vessel. 
The stenosis degree is parameterized by two geometric parameters: width, $w \in \{0.4,0.6\}cm$, and depth, $d \in \{0.1,0.3,0.5,0.7\}cm$. 


As for P-inflow, stability is almost completely preserved across all geometric configurations when the original network (trained for the non stenotic case) is used to estimate $\widehat{\alpha}$; see Table \ref{tab:stability}.
Fig. \ref{fig:errors_network}-C shows that the median of the relative error distributions is lower than in Fig. \ref{fig:errors_network}-A, indicating an overall increase in the accuracy for the LC solutions.

These results demonstrate that the predicted stabilization parameter generalizes well with respect to geometric variation of the computational domain.
Moreover, they suggest that introducing an obstruction to the flow enhances the performance of the proposed partitioned algorithm, leading to higher accuracy than that obtained with the baseline problem (straight rectangle) same value of \al.
This could be due to a beneficial effect of the flow restriction that somehow partially compensates for the added-mass effect.

\begin{figure}
  \centering
  \includegraphics[width=0.99\textwidth]{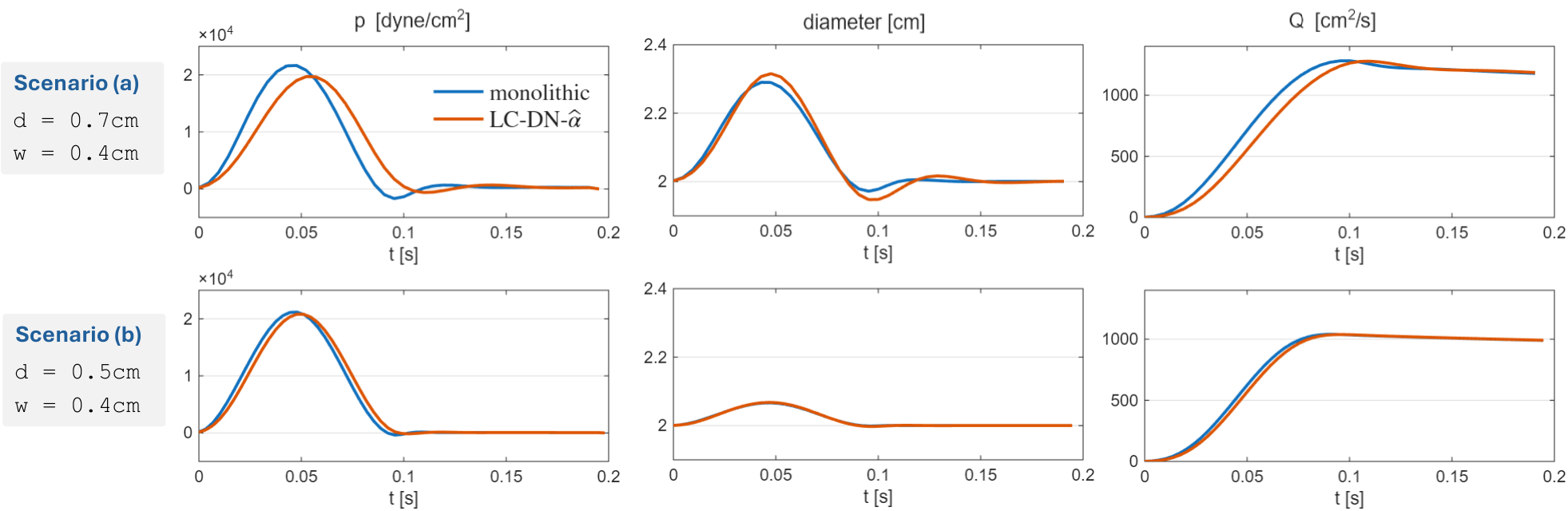}
  \caption{Comparison between the evolution over time of the LC-DN-$\widehat{\alpha}$ and the monolithic solutions, orange and blue lines, respectively, extracted at the mid-section of the domain ($x=3cm$) for two stenotic scenarios. 
  The physical and numerical characteristics of the reported scenarios are detailed in Fig.\ref{fig:curves_cons}.
  The gray panels on the left show the geometric parameters that define the degree of stenosis. From the left: the first column shows the fluid pressure along the y-direction, the second the diameter of the fluid domain, and the third the fluid flow. Robustness with respect to changes in the BP: P-stenosis.}
  \label{fig:stenosis_curves}
\end{figure}

\begin{figure}
  \centering
  \includegraphics[width=0.99\textwidth]{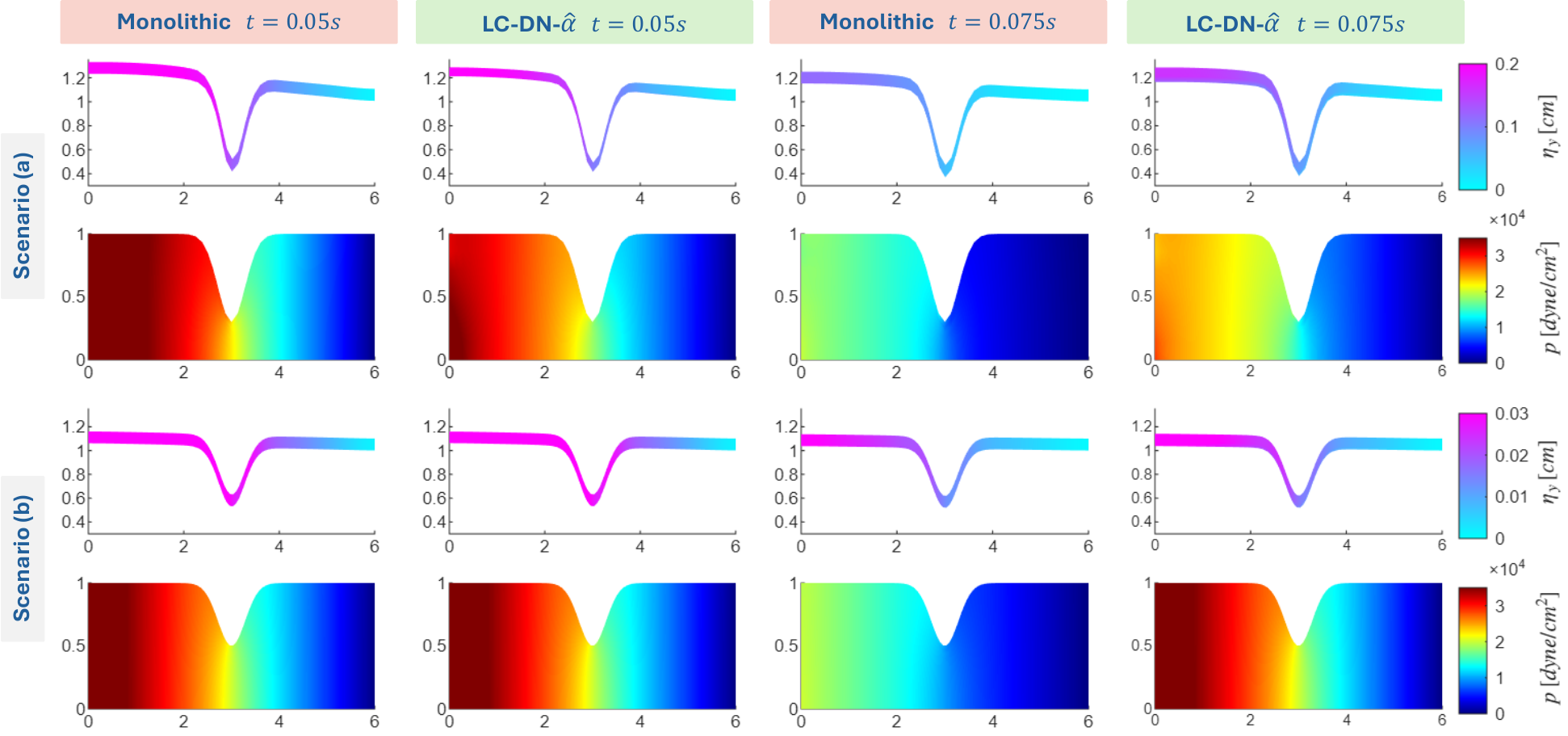}
  \caption{Plots of the LC-DN-$\widehat{\alpha}$ and the monolithic solutions extracted at two time instants, for the same two FSI scenarios as in \ref{fig:stenosis_curves}. 
  For each scenario, the structure displacement along x and the fluid pressure (top and bottom rows, respectively) are reported. Robustness with respect to changes in the BP: P-stenosis.}
  \label{fig:stenosis_plots}
\end{figure}

\begin{figure}
  \centering
  \includegraphics[width=0.6\textwidth]{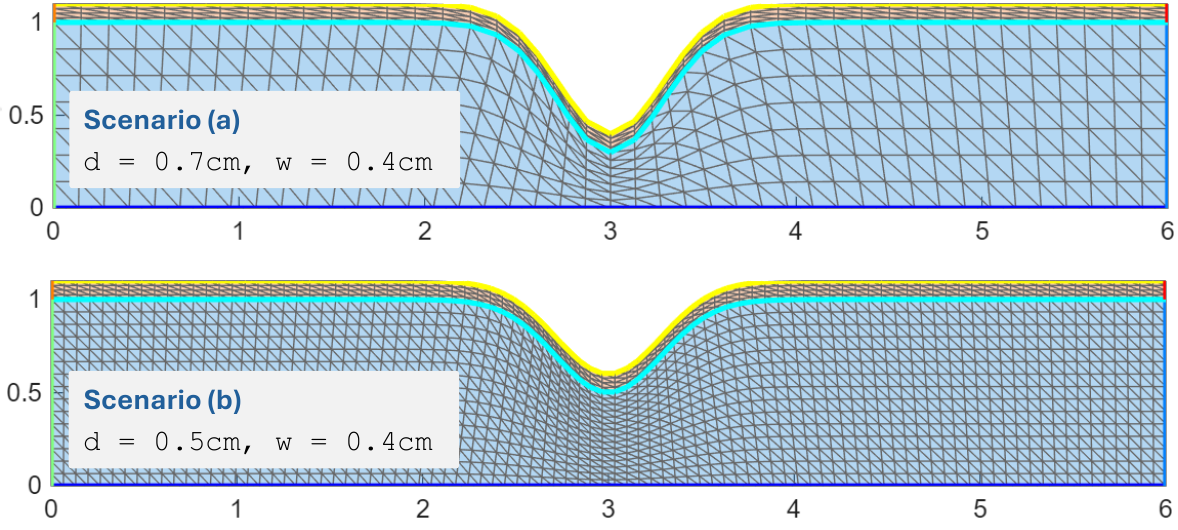}
  \caption{Examples of modified geometries for the scenarios (a) and (b) (top and bottom figures, respectively). Robustness with respect to changes in the BP: P-stenosis.}
  \label{fig:stenosis_examples}
\end{figure}

Finally, in Figures \ref{fig:stenosis_curves} and \ref{fig:stenosis_plots}, we report the computed solution at two time instants for scenarios (a) and (b). The characteristics of these scenarios are detailed in Fig. \ref{fig:curves_cons}, while their geometric configurations are displayed in Fig. \ref{fig:stenosis_examples}.
These results confirm the above consideration regarding accuracy with respect to the reference solution.

\begin{table}[ht]
\centering
\caption{Percentage of stable simulations for each combination of stenosis depth and width. Robustness with respect to changes in the BP: P-stenosis.}
\label{tab:stability}
\begin{tabular}{|c|c|c|}
\toprule
\diagbox[width=7em]{$d [cm]$}{$w [cm]$} & 0.4 & 0.6 \\
\midrule
0.1 & 100\% &  100\% \\
0.3 & 100\% &  100\% \\
0.5 & 97.8\% & 97.8\% \\
0.7 & 95.6\% & 95.6\% \\
\bottomrule
\end{tabular}
\end{table}

\section{Conclusions}
In this work, we propose a framework for the automatic identification of the stabilization parameter of the LC-DN-\al scheme that exploits the capabilities of SciML. 

The analysis performed in Section \ref{subsec:features} highlights nonlinear relationships among the physical, numerical, and algebraic features of the FSI problems, as well as between these features and the stabilization parameter $\alpha$, thus justifying the need for a nonlinear model to express the relationship between the problem features and \al. 
Based on this analysis, we also concluded that the time step $\Delta t$, the structure Young's modulus $E$, and the structure density $\rho_s$ are sufficient as input features for the neural network. 
The resulting model was shown to generalize well across the considered test cases (Section \ref{subsec:consistency}) and also to provide robust parameter estimates under moderate variations of the problem setup (Section \ref{subsec:robustness}), specifically with respect to the loading conditions (P-inflow, Fig. \ref{fig:fxin}) and computational domain (P-stenosis, Fig. \ref{fig:errors_network}-C).
This suggest that it is not necessary to re-train the network for these two kind of modifications.
Degradation in the performance of the network predictions is,  however, observed with respect to significant changes in the outlet boundary conditions (P-outflow, Fig. \ref{fig:errors_network}-B).

Nevertheless, the overall robustness demonstrated by the proposed approach is particularly encouraging for hemodynamic applications.
Indeed, as simulations become more complex and patient-specific, variations in geometry and inlet boundary conditions are common from one patient to another. 
In this context, the model's ability to provide reliable parameter estimates represents an important achievement toward reducing the computational demand, with the broader goal of translating patient-specific computational hemodynamics into clinical practice.

The proposed methodology has been intentionally developed and validated for simplified FSI problems involving Stokes flows and linear-elastic structures in 2D domains.
Indeed, at this stage, this study is intended to provide a foundation for extending it to more complex and realistic scenarios. 
Future developments must shift toward Navier–Stokes equations and nonlinear constitutive models for the structure.
Such extension is not straightforward, as the relationship between the FSI problem features and the optimal stabilization parameters may become significantly more complex than that observed in the present study.
Consequently, the current set of input features used for the network's training ($\Delta t, E, \rho_s$) may no longer be sufficient to characterize the problem's variability, necessitating the inclusion of additional descriptors of the coupled system.
Furthermore, the extension to 3D configurations represents a necessary step toward real-world applications and will provide a critical assessment of the general validity of the proposed approach.

\section*{Acknowledgments}
\noindent FR, SP, PFA, CV are members of INdAM GNCS. 

\noindent FR is partially supported by the INdAM GNCS project 2026 - CUP E53C25002010001.
FR has been funded by the European Union-Next Generation EU, Mission 4, Component 1, CUP: D53D23018770001, research project MIUR PRIN22-PNRR n.P20223KSS2.

\noindent CV has been partially supported by i) MIUR PRIN22-PNRR n.P20223KSS2; ii) Italian Ministry of Health, PNC PROGETTO HUB LIFE SCIENCE - DIAGNOSTICA AVANZATA (HLSDA) \lq\lq INNOVA\rq\rq, PNCE3 - 2022-23683266-CUP: D43C22004930001; iii) Research project MIUR PRIN22 n.2022L3JC5T; iv) Italian Ministry of Health, project CALHUBRIA, Codice Progetto: T4-AN-09, CUP: F63C22000530001. 

\noindent PFA ans SP are partially funded by the European Union (ERC SyG, NEMESIS, project number 101115663).
Views and opinions expressed are, however, those of the authors only and do not necessarily reflect those of the European Union or the European Research Council Executive Agency. Neither the European Union nor the granting authority can be held responsible for them. 

\noindent The present research is part of the activities of the Dipartimento di Eccellenza 2023-2027 grant, funded by MUR.

\noindent We acknowledge the CINECA award under the ISCRA initiative, for the availability of high-performance computing resources and support.

\bibliographystyle{cas-model2-names}

\bibliography{Renzi-et-al-CMAME2026_references}

\end{document}